\documentclass[aap,MSNbibl,citesort,seceqn,dvips]{arximspdf}
\usepackage{graphicx}

\doi{10.1214/11-AAP809} %kopijuoti is PTS
\volume{22}
\issue{4}
\pubyear{2012}
\firstpage{1693}
\lastpage{1711}

\makeatletter

\newtheorem{teor}{Theorem}
\newtheorem{prop}[teor]{Proposition}
\newtheorem{lem}[teor]{Lemma}
\newtheorem{cor}[teor]{Corollary}

\newproclaim{Def}[teor]{Definition}
\newproclaim{cond}[teor]{Condition}
\newproclaim{conj}[teor]{Conjecture}
\newproclaim{Exa}[teor]{Example}
\newproclaim{rem}[teor]{Remark}
\newcommand{\R}{\mathbb{R}}
\newcommand{\N}{\mathbb{N}}
\newcommand{\PP}{\mathbb{P}}
\newcommand{\E}{{\mathbb{E}}}
\newcommand{\defi}{\equiv} % {\stackrel{\text{def}}{=}}
\newcommand{\be}{\beta}
\newcommand{\de}{\delta}
\newcommand{\dd}{\mathrm{d}}
\newcommand{\ee}{\mathrm{e}}
\newcommand{\w}{\omega}
\newcommand{\vare}{\varepsilon}
\newcommand{\eqref}[1]{(\ref{#1})}
\makeatother

\begin{document}
\begin{frontmatter}

\title{Poissonian statistics in the extremal process of branching Brownian motion\thanksref{T0}}
\runtitle{Poissonian statistics in branching Brownian motion}
\thankstext{T0}{Part of this work has been carried out during the
\textit{Junior Trimester Program Stochastics} at the
Hausdorff Center in Bonn; hospitality and financial support are
gratefully acknowledged.}

\begin{aug}
\author[A]{\fnms{Louis-Pierre} \snm{Arguin}\corref{}\ead[label=e1]{arguinlp@dms.umontreal.ca}\thanksref{t1}},
\author[B]{\fnms{Anton} \snm{Bovier}\ead[label=e2]{bovier@uni-bonn.de}\thanksref{t2}}
\and
\author[B]{\fnms{Nicola} \snm{Kistler}\ead[label=e3]{nkistler@uni-bonn.de}\thanksref{t3}}
\runauthor{L.-P. Arguin, A. Bovier and N. Kistler}
\thankstext{t1}{Supported by NSF Grant DMS-06-04869.}
\thankstext{t2}{Supported in part through the German Research Council
in the SFB 611 and
the Hausdorff Center for Mathematics.}
\thankstext{t3}{Supported in part by the Hausdorff Center for Mathematics.}
\affiliation{Universit\'e de Montr\'eal, Rheinische
Friedrich-Wilhelms-Universit\"at
Bonn and Rheinische Friedrich-Wilhelms-Universit\"at Bonn}
\address[A]{L.-P. Arguin\\
D\'epartement de Math\'ematiques \\
\quad et Statistique\\
Universit\'e de Montr\'eal\\
C.P. 6128, succ. Centre-ville\\
Montr\'eal, Qu\'ebec\\
Canada H3C 3J7\\
\printead{e1}} %adresu isvedimo komanda gale!
\address[B]{A. Bovier\\
N. Kistler\\
Institut f\"ur Angewandte Mathematik\\
Rheinische Friedrich-Wilhelms-Universit\"at Bonn\\
Endenicher Allee 60\\
53115 Bonn\\
Germany\\
\printead{e2}\\
\phantom{E-mail: }\printead*{e3}}
\end{aug}

% HISTORY:
\received{\smonth{11} \syear{2010}}
\revised{\smonth{8} \syear{2011}}

% ABSTRACT
%
\begin{abstract}
As a first step toward a characterization of the limiting extremal
process of branching Brownian motion, we proved in a recent work
[\textit{Comm. Pure Appl. Math.} \textbf{64} (2011) 1647--1676]
that, in the limit of large time $t$, extremal particles descend
with overwhelming probability from ancestors having split either within
a distance of order~1 from time $0$, or within a distance of order~1
from time $t$. The result suggests that the extremal process of
branching Brownian motion is a randomly shifted cluster point process.
Here we put part of this picture on rigorous ground: we prove that the
point process obtained by retaining only those extremal particles which
are also maximal inside the clusters converges in the limit of large
$t$ to a random shift of a~Poisson point process with exponential
density. The last section discusses the \textit{Tidal Wave Conjecture} by
Lalley and Sellke
[\textit{Ann. Probab.} \textbf{15} (1987) 1052--1061] on the full
limiting extremal
process and its relation to the work of Chauvin and Rouault
[\textit{Math. Nachr.} \textbf{149} (1990) 41--59] on branching
Brownian motion with atypical displacement.
\end{abstract}

% KEYWORDS
%

\end{frontmatter}

%s1 ###
\section{Introduction}\label{sec1}
Branching Brownian motion (BBM) is a continuous-time Markov branching process
which plays an important role in the theory of partial differential
equations~\cite{aronsonweinberger,aronsonweinbergertwo,mckean},
in particle physics~\cite{munierpeschanski}
and in the theory of disordered systems
\cite{BovierKurkovaII,derridaspohn}. It is also widely used in
biology to model the genealogies of evolving populations, the spread of
advantageous genes, etc.,~\cite{fisher,kessleretal}. It is
constructed as follows.

Start with a standard Brownian motion (BM) (we will often
refer to Brownian motions as ``particles''), $x(t)$, starting at $0$.
After an exponential random time, $T$, of mean 1, the BM splits into $k$
independent BMs, independent of $x$ and $T$, with probability $p_k$,
where $\sum_{k=1}^\infty p_k = 1$, $\sum_{k=1}^\infty k p_k = 2$ and
$K \defi\sum_{k} k(k-1) p_k < \infty$. Each of these processes
continues in the same way as first BM. Thus, after time $t>0$, there
will be $n(t)$ BMs located at
$x_1(t), \ldots, x_{n(t)}(t)$, with $n(t)$ being the random number of offspring
generated up to that time [note that $\E n(t)=e^t$].

An interesting link between BBM and partial differential equations
was observed by McKean~\cite{mckean}: denote by
\begin{equation}\label{bbmrepr}
u(t, x) \defi\PP\Bigl[ \max_{1\leq k \leq n(t)} x_k(t) \leq x \Bigr]
\end{equation}
the law of the maximal displacement. Then, a renewal argument shows that
$u(t,x)$ solves
the Kolmogorov--Petrovsky--Piscounov or Fisher [F-KPP] equation,
\begin{eqnarray}\label{kppequation}
 u_t &=& {1\over2} u_{xx} + \sum_{k=1}^\infty p_k u^k -u,
 \nonumber
 \\[-8pt]
 \\[-8pt]
 \nonumber
 u(0, x)&=&
\cases{
1, &\quad$\mbox{if } x\geq0,$\vspace*{2pt}\cr
0, &\quad$\mbox{if } x < 0.$}
\end{eqnarray}
The F-KPP equation admits traveling waves: there exists a unique
solution satisfying
\begin{equation}\label{travellingone}
u\bigl(t, m(t)+ x \bigr) \to\omega(x)\qquad \mbox{uniformly in } x
\mbox{ as } t\to\infty,
\end{equation}
with the centering term, the \textit{front} of the wave, given by
\begin{equation}\label{centeringkpp}
m(t) = \sqrt{2} t - {3 \over2 \sqrt{2}} \log t, %+ b(t), b(t) =
%O(1) \mbox{for} t\to\infty,
\end{equation}
and $\w(x)$ the unique (up to translation) distribution function which
solves the ordinary differential equation
\begin{equation}\label{wavepde}
\tfrac{1}{2} \omega_{xx} + \sqrt{2} \omega_x + \omega^2 - \omega= 0.
\end{equation}
The leading order of the front has been established by Kolmogorov,
Petrovsky and Piscounov~\cite{kpp}, whereas the logarithmic corrections
have been obtained
by Bramson~\cite{bramson}, using the probabilistic representation given
above.

The limiting law of the maximal displacement has been studied intensely.
Let
%e1 ###
%
\begin{equation}
\label{eqnderiv}
Z(t) \defi\sum_{k=1}^{n(t)} \bigl( \sqrt{2}t -x_k(t) \bigr) \exp-\sqrt
{2}\bigl( \sqrt{2}t -x_k(t) \bigr)
\end{equation}
denote the so-called \textit{derivative martingale}.
Lalley and Sellke~\cite{lalleysellke} proved that~$Z(t)$ converges almost
surely to a strictly positive
random variable, $Z$, and established the integral representation
\begin{equation}\label{gumbellike}
\omega(x) = \E\bigl[ \ee^{- C Z \ee^{-\sqrt{2}x} }\bigr],
\end{equation}
with $C>0$ a constant. Thus the law of the maximum of BBM is a
\textit{random shift} of the Gumbel distribution. Moreover, it is
known that
\begin{equation}\label{totheright}
1-\omega(x) \sim x \ee^{-\sqrt{2} x},\qquad x\to+\infty,
\end{equation}
where $\sim$ means that the ratio of the terms converges to a positive
constant; see, for example, Bramson~\cite{bramson} and Harris~\cite{harris}.
(There is emerging evidence that right-tails such as \eqref{totheright},
manifestly different from those of the Gumbel, play an important
role in a number of different fields, e.g., in models on spin glasses
with logarithmic correlated potentials by Carpentier and Le Doussal~\cite{carpentierledoussal}, and Fyodorov and Bouchaud \cite
{bouchaudfyodorov}.)

Contrary to the maximal displacement, very little is known on the
full statistics of the extremal configurations (first-, second-,
third-, etc., largest) in BBM.
Such statistics are completely encoded in the extremal process, which
is the random point measure associated to the collection
of points shifted by the expectation of their maximum, that is,
the point process
%e2 ###
%
\begin{equation} \label{extremalprocess}
\Xi(t) \defi\sum_{i=1}^{n(t)} \de_{\overline{x_i}(t)} ,\qquad
\overline{x_i}(t) \defi x_i(t)-m(t).
\end{equation}
The key issue of interest is to characterize the limit of this process,
as $t\uparrow\infty$. It can be shown that the limit of the point process
exists using Bramson's analysis~\cite{bramsonmonograph} on the
convergence of solutions of the KPP equations with appropriate initial
conditions~\cite{brunet,kabluchko}.

For given realization of the branching, the positions
$\{x_i(t)\}_{i\leq n(t)}$
form a~Gaussian process indexed by $i\in\{1, \ldots, n(t)\} \defi
\Sigma_t$
with correlations
given by the \textit{genealogical distance}
\begin{equation}
Q_{ij}(t) \defi\sup\{s\leq t\dvtx  x_i(s)=x_j(s) \} \label{genealogical}
\end{equation}
(the time to first branching of the common ancestor). The information
about the correlation structure of any subsets of particles in BBM
is encoded in their genealogical distance. This applies, in particular,
to the subset of extremal particles, for which the following result was proved
in~\cite{abk}: with probability tending to 1, branching can happen
only at ``very early times,'' smaller than $r_d$ with $r_d = O(1)$ as
$t\to\infty$, or at times ``very close'' to the age of the system,
namely greater than
$t-r_g$ for
$r_g = O(1)$ as $t \to\infty$. (The reason for this notation,
in particular the use of the subscripts, will be explained below.)
More precisely, denoting by $\Sigma_t(D) \defi\{ i\in\Sigma_t\dvtx
\overline{x_i}(t) \in D \}$ the set of particles in the subset
$m(t)+D$, we have:
\begin{teor}[\cite{abk}]
\label{genealogyabktheorem}
For any compact $D\subset\R$,
\begin{equation}\label{genealogyabkeqn}
\lim_{r_d, r_g \to\infty} \sup_{t > 3 \max\{r_d, r_g \}} \PP[
\exists i,j \in\Sigma_t(D)\dvtx  Q_{ij}(t) \in(r_d, t-r_g) ] = 0.
\end{equation}
\end{teor}

\begin{figure}

\includegraphics{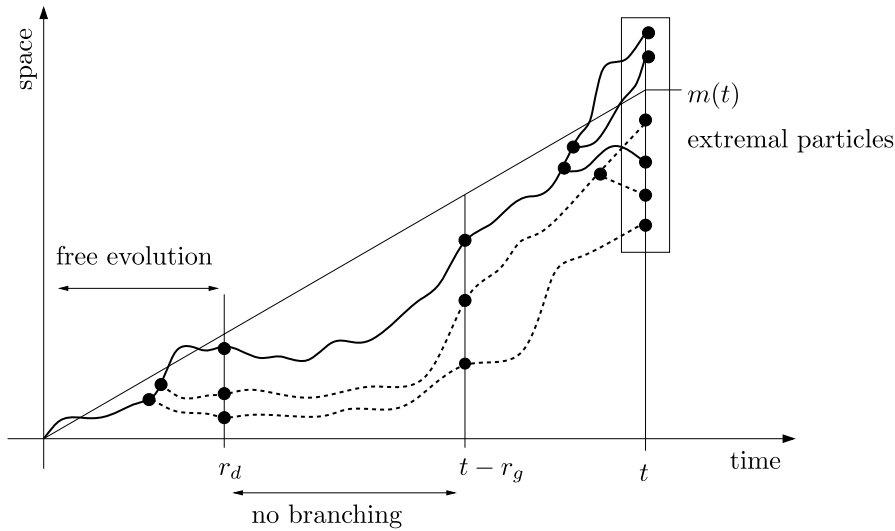}

\caption{Genealogies of extremal particles.}\label{fig1}
\end{figure}

Figure~\ref{fig1} presents a graphical representation of the
genealogies of
extremal particles of BBM.

%%f1 ###
%%
%%
%%
%%
%%
%%
%{10}{12.0}{\rmdefault}{\mddefault}{\updefault}{
%%%
%}}}}
%{10}{12.0}{\rmdefault}{\mddefault}{\updefault}{\color[rgb]{0,0,0}$t$}%
%}}}}
%{10}{12.0}{\rmdefault}{\mddefault}{\updefault}{\color
%[rgb]{0,0,0}$t-r_g$}%
%}}}}
%{10}{12.0}{\rmdefault}{\mddefault}{\updefault}{\color[rgb]{0,0,0}free
%evolution}%
%}}}}
%{10}{12.0}{\rmdefault}{\mddefault}{\updefault}{
%}}}}
%{10}{12.0}{\rmdefault}{\mddefault}{\updefault}{\color[rgb]{0,0,0}no
%branching}%
%}}}}
%{10}{12.0}{\rmdefault}{\mddefault}{\updefault}{\color
%[rgb]{0,0,0}extremal particles}%
%}}}}
%{{\SetFigFontNFSS{10}{12.0}{\rmdefault}{\mddefault}{\updefault}{\color
%[rgb]{0,0,0}space}%
%}}}}}
%{10}{12.0}{\rmdefault}{\mddefault}{\updefault}{\color[rgb]{0,0,0}time}%
%}}}}
%%
%
%%
%

% \input{picture.pstex_t}

Theorem~\ref{genealogyabktheorem} gives insight into the limiting extremal
process of BBM. In fact, it suggests the following picture, which holds with
overwhelming probability in the limit when first $t\uparrow\infty$, and
$r_d,r_g\rightarrow\infty$ after that.

First, ancestries in the interval $[0, r_d]$ cannot be ruled out:
this regime generates the derivative martingale appearing in the work of
Lalley and Sellke~\cite{lalleysellke}. Moreover, since the ancestors
of the
extremal particles
evolved independently for most of the time (namely
in the interval $[r_d,t-r_g]$), the extremal process must exhibit a
structure similar to a Poisson process. Finally, since ancestors over the
period
$[t-r_g, t]$ also occur, it is natural
to conjecture that \textit{small grapes} of length at most $r_g = O(1)$,
that is,
clusters of particles with very recent common ancestor,
appear at the end of the time-interval.
(According to this picture, the subscript in $r_d$ refers to \textit{
derivative martingale}, while that in $r_g$ stands for \textit{grape}.)

It is the purpose of this work to make part of this picture rigorous.
In Section~\ref{mainresults} we present our main result,
which is proved in Section~\ref{proofs}. In Section~\ref{secconvection},
we introduce a cluster point process, which we conjecture to correspond in
the limit
to the extremal process of BBM. We also discuss the cluster point
process in
relation to the work of Chauvin and Rouault~\cite{chauvinrouault} on BBM
conditioned to perform unusually large displacements, and in relation
to the
\textit{Tidal Wave Conjecture} of Lalley and Sellke~\cite{lalleysellke}.
Detailed properties of this cluster point process will be the subject
of a
subsequent paper~\cite{abkthree}.

%s2 ###
\section{Main results} \label{mainresults}

Despite the rather clear image described above, a frontal attack on the
extremal process appears to be difficult.
This is in particular due to the fact that one has to take into account the
self-similarity of BBM which is first and foremost detectable in the
small clusters, an issue which remains rather elusive
(see Section~\ref{secconvection} for more on this).
On the other hand, the picture naturally suggests the existence of an
underlying point process obtained from the extremal particles
by a thinning procedure, which we describe next.

Assume that the positions of particles at time $t$ are ordered in
decreasing order:
\begin{equation}
\overline{x}_1(t)\geq\overline{x}_2(t)\geq\cdots\geq\overline
{x}_{n(t)}(t).
\end{equation}
The inequalities will in fact be strict for almost all realizations of BBM
for any deterministic time $t$. Define also
\begin{equation}
\overline{Q}(t) =\{\overline{Q}_{ij}(t)\}_{i,j\leq n(t)}\defi\{
t^{-1}Q_{ij}(t)\}_{i,j\leq n(t)} .
\end{equation}
The pair $(\Xi(t), \overline{Q}(t))$ admits the following natural
\textit{thinning}.
Since the matrix $\overline{Q}(t)$ is constructed from the branching of
the BBM, the relation
$\overline{Q}_{ij}(t)\geq q$ is transitive for any $q\geq0$:
%e3 ###
%
\begin{equation}
\label{eqnultra}
\overline{Q}_{ij}(t)\geq q \quad\mbox{and} \quad \overline{Q}_{jk}(t)\geq q
\quad\Longrightarrow\quad\overline{Q}_{ik}(t)\geq q.
\end{equation}
In particular, for any $q>0$, this relation defines an equivalence relation
on the set $\{1,\ldots,n(t)\}$. The corresponding equivalence classes are just
the
particles at time~$t$ that had a common ancestor at a time later than $tq$.
We want to select a representative of each class, namely the maximal
particle within each class, and then consider the point process of these
representatives.
For any $q>0$, the \mbox{\textit{$q$-thinning}} of the process $(\Xi(t),
\overline
{Q}(t))$, denoted by $\Xi^{(q)}(t)$, is defined recursively as follows:
%e4 ###
%
\begin{eqnarray}
\label{eqnthin}
i_1&=&1;
\nonumber
\\[-8pt]
\\[-8pt]
\nonumber
i_k&=&\min\{ j>i_{k-1}\dvtx  \overline{Q}_{i_lj}(t)<q,  \forall l\leq k-1\};
\end{eqnarray}
and
%e5 ###
%
\begin{equation}
\Xi^{(q)}(t) \defi\bigl(\Xi^{(q)}_{k}(t), k\in\N\bigr) \defi\bigl(\overline
{x}_{i_k}(t), k\in\N\bigr),
\end{equation}
where it is understood that $\Xi^{(q)}_{k}(t)=0$ when an index $i_k$ in
$\{1,\ldots,n(t)\}$ satisfying $\min\{ j>i_{k-1}\dvtx  \overline{Q}_{i_lj}<q
\ \forall l\leq k-1\}$ can no longer be found.
The procedure selects the maximal position in each equivalence class
defined from the relation \mbox{$\overline{Q}_{ij}(t)\geq q$}.
In addition, it is easily checked that the thinning map,
%e6 ###
%
\begin{equation}
\label{eqnthinmap}
(\Xi(t),\overline{Q}(t))\mapsto\Xi^{(q)}(t),
\end{equation}
considered at the level of realizations, is a continuous function on the
space of pairs $(X,Q)$, where $X$ is
a sequence of ordered positions and $Q$ is a~symmetric matrix with entries
in $[0,1]$, satisfying \eqref{eqnultra}
(when this space is equipped with the product topology in each
coordinate of $X$ and $Q$).

\begin{figure}

\includegraphics{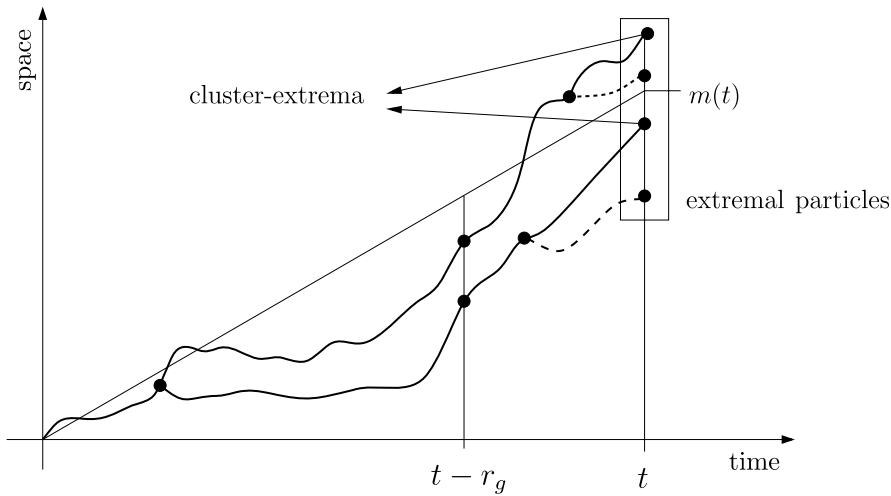}

\caption{Cluster-extrema.}\label{fig2}
\end{figure}

The thinning map can also be applied to $t$-dependent values of $q$.
For example, take $q=q(t)=1-r_g/t$, where $r_g$ is fixed $t$.
In this case, the thinning effectively retains those particles which
are extremal within the class defined by a ``very recent'' common
ancestor, which we refer to as
\textit{cluster-extrema}. Figure~\ref{fig2} presents a graphical
representation of the set of such particles.

%f2 ###

%%
%%
%%
%%
%%
%%
%{{\SetFigFontNFSS{10}{12.0}{\rmdefault}{\mddefault}{\updefault
%}{\color
%[rgb]{0,0,0}space}%
%}}}}}
%{10}{12.0}{\rmdefault}{\mddefault}{\updefault}{\color
%[rgb]{0,0,0}time}%
%}}}}
%{12}{14.4}{\rmdefault}{\mddefault}{\updefault}{\color[rgb]{0,0,0}$m(t)$}
%%%
%}}}}
%{12}{14.4}{\rmdefault}{\mddefault}{\updefault}{\color
%[rgb]{0,0,0}$t-r_g$}%
%}}}}
%{12}{14.4}{\rmdefault}{\mddefault}{\updefault}{\color
%[rgb]{0,0,0}$t$}%
%}}}}
%{12}{14.4}{\rmdefault}{\mddefault}{\updefault}{\color
%[rgb]{0,0,0}extremal particles}%
%}}}}
%{12}{14.4}{\rmdefault}{\mddefault}{\updefault}{\color
%[rgb]{0,0,0}cluster-extrema}%
%}}}}
%%
%
%%
%

Our main result states that all such thinned processes converge to the
same randomly shifted Poisson Point Process (PPP for short) with
exponential density.
\begin{teor}
\label{convergencethinned}
For any $0<q<1$,
the processes $\Xi^{(q)}(t)$ converge in law to the same limit, $\Xi
^0$. Also,
\begin{equation}
\lim_{r_g \to\infty} \lim_{t\to\infty} \Xi^{(1-r_g/t)}(t)=\Xi^0.
\end{equation}
Moreover, conditionally on $Z$, the limit of the derivative martingale
\eqref{eqnderiv},
\begin{equation}
\Xi^0=\operatorname{PPP}\bigl( C\cdot Z \cdot\sqrt{2}\ee^{-\sqrt{2} x} \,\dd
x\bigr) ,
\end{equation}
where $C>0$ is the constant appearing in \eqref{totheright}.
\end{teor}

The point process $\Xi^{0}$ has a fundamental connection with the limiting
extremal process of BBM.
To see this, suppose for simplicity that the processes $(\Xi(t),
\overline{Q}(t))$ induced by the law of BBM
converge, as $t\uparrow\infty$, to a process, $(\Xi, \overline{Q})$.
(The laws of these processes are in fact tight because the law of $\Xi
(t)$ is
itself tight; see, e.g., Corollary 2.3 in~\cite{abk},
and that $\overline{Q}_{ij}(t)\in[0,1]$ for any~$i,j$.
Convergence would evidently follow from a complete characterization of the
extremal process.)
It follows from Theorem~\ref{genealogyabktheorem} that $\overline{Q}_{ij}$
is either $0$ or $1$.
This suggests:
\begin{enumerate}[(1)]
\item[(1)] to define a cluster of particles as the maximal set of particles such
that $\overline{Q}_{ij}=1$ for all $i,j$ in the set;
\item[(2)] to look at the process of the maxima of each cluster, denoted by, say,
$\tilde{\Xi}^0$, defined as in \eqref{eqnthin}, but where
$i_k=\min
\{ j>i_{k-1}\dvtx  \overline{Q}_{i_lj}=0\  \forall l\leq k-1\}$.
\end{enumerate}
We claim that $\tilde{\Xi}^0$ is in fact the limit $\Xi^{0}$ of $\Xi
^{(q)}(t)$ in Theorem~\ref{convergencethinned}.
Indeed, in view of the continuity of the thinning map \eqref{eqnthinmap}, $\Xi^{(q)}(t)$
converges to the \mbox{$q$-thinned} process, $\Xi^{(q)}$,
constructed from $(\Xi, \overline{Q})$ for all $q$.
But, for any $0<q<1$, the \mbox{$q$-thinned} processes, $\Xi^{(q)}$,
constructed from $(\Xi, \overline{Q})$ using \eqref{eqnthin} are
equal trivially to $\tilde{\Xi}^0$, since $\overline{Q}_{ij}$ is either
$0$ or $1$. The claim then follows from Theorem~\ref{convergencethinned}.
The point process describing the particles at the frontier of BBM in
the limit of large times is thus formed by two ``types'' of particles: those
coming from the randomly shifted PPP with exponential density, the
cluster-extrema; and the second type of particles, those forming the clusters.
Clearly, particles in the same cluster always lie on the left of the
corresponding Poissonian particles, by the very definition of the
cluster-extrema.
It remains an open question to characterize the law of the clusters
(see Section~\ref{secconvection} for some conjectures).

We remark that, since $r_g = O(1)$ as $t\to\infty$, the thinned
process $\Xi^{(1-r_g/t)}_t$ is obtained from the extremal one by removing
only a small number of particles, those which have genealogical
distance smaller than $t-r_g$ from the maximum in their class. It is
rather surprising at first sight (but not quite when seen under the
light of Theorem~\ref{genealogyabktheorem}) that such a point process
converges, despite the high correlations among the branching Brownian
particles, to a~PPP with exponential density.

Theorem~\ref{convergencethinned} also provides insights into a result
by Bovier and Kurkova~\cite{BovierKurkovaII}, who addressed the weak
limit of the Gibbs measure
of BBM, the random probability measure on $\Sigma_t$ attaching weights
\begin{equation}
\mathcal G_{\beta, t}(k) \defi\frac{\exp(\beta x_k(t))}{\mathcal
{Z}_t(\beta)},\qquad \mathcal{Z}_t(\beta) \defi\sum_{j\in\Sigma_t}
\exp(\beta x_j(t)),
\end{equation}
where $\beta>0$ is the inverse of temperature. To see this, let us
first recall the following.

Consider the random set $(\xi_i, i \in\N)$ where the $\xi$'s are
generated according to a PPP with density $C Z \sqrt{2}\ee^{-\sqrt{2}x}
\,\dd x$ on the real axis, $C$ and $Z$ as in Theorem~\ref{convergencethinned}.
Construct then a new random set $(\rho_i, i\in
\N)$ where $\rho_i \defi\exp(\be\xi_i)$. For $\beta> \sqrt{2}$, it is
easily seen that $\mathcal N(\rho) \defi\sum_j \rho_j < \infty$ almost
surely, in which case the normalization $\hat{\rho}_i \defi\rho
_i/\mathcal N(\rho)$ is well defined, and the law of the normalized
collection $(\hat{\rho}_i, i \in\N)$ is the Poisson--Dirichlet
distribution with parameter $m(\beta) = \sqrt{2}/\beta$, which we shall
denote by $\operatorname{PD}(m(\beta))$.

In a somewhat indirect way (by means of the so-called Ghirlanda--Guerra
identities, which avoid the need to first identify the limiting
extremal process) Bovier and Kurkova proved that, in the low
temperature regime $\beta> \sqrt{2}$, the Gibbs measure $\mathcal
G_{\beta, t}$ converges, in the limit of large times, to the
$\operatorname{PD}(m(\beta
))$; together with our Theorem~\ref{convergencethinned}, this
naturally suggests that the Gibbs measure of BBM is concentrated, in
fact, on the \textit{cluster-extrema}.

Finally, Theorem~\ref{convergencethinned} sheds light on a property of
the extremal process of BBM that was conjectured by Brunet and Derrida
\cite{derridabrunet}.
They suggested that the statistics of the leading particles are
invariant under superposition in the sense that the extremal process of
two independent branching Brownian motions has the same law, up to a
random shift, as the extremal process
of a single one. This property at the level of the entire process is
likely to involve specific features of the laws of the individual
clusters. On the other hand,
at the level of the thinned process, it is a straightforward
consequence of Theorem~\ref{convergencethinned}, since the law is
Poisson with exponential density.
\begin{cor}
\label{corsuperposition}
Let $\Xi(t)$ and $\Xi'(t)$ be the extremal processes \eqref{extremalprocess} of two independent branching Brownian motions.
Denote by $Z$ and $Z'$ the pointwise limit of their respective
derivative martingale.
Then, for any $0<q<1$, the law of the $q$-thinning of $\Xi(t)+\Xi'(t)$
conditionally on $Z$ and $Z'$ converges to
\begin{equation}
\operatorname{PPP}\bigl( C\cdot(Z+Z')\cdot\sqrt{2}\ee^{-\sqrt{2} x} \,\dd
x\bigr) .
\end{equation}
In particular, the thinned process of $\Xi(t)+\Xi'(t)$ has the same law
in the limit as the thinned process $\Xi^0$ of a single
branching Brownian motion, up to a random shift.
\end{cor}

As mentioned before, Theorem~\ref{convergencethinned} is a natural
consequence of Theorem~\ref{genealogyabktheorem}.
The main ingredient is the following lemma, which allows to compare thinning
processes on a set of large probability.
We use the notation
\begin{equation}
\Xi^{(q)}(t)|_{(y,\infty)}\defi\bigl\{\Xi^{(q)}_i(t)\dvtx  \Xi^{(q)}_i(t)>y\bigr\}
 .
\end{equation}

\begin{lem}
\label{comparison}
For any $y\in\R$ and any $\vare>0$, there exists $r_0=r_0(y,\vare)$
such that for $r_d,r_g>r_0$ and $t>3\max\{r_g,r_d\}$, on a set of
probability $1-\vare$,
\begin{equation}
\Xi^{(q)}(t)|_{(y,\infty)}= \Xi^{(r_d/t)}_t|_{(y,\infty)},
\end{equation}
for any $\frac{r_d}{t}<q<1-\frac{r_g}{t}$.
\end{lem}

Theorem~\ref{convergencethinned} is then proved by a standard Poisson
convergence argument which exploits the weak correlations between
the cluster-extrema in classes of the $\frac{r_d}{t}$-thinning.\vspace*{-3pt}

\begin{prop} \label{simplerpp}
With $C>0$ and $Z$ the limiting derivative martingale, conditionally on $Z$,
\begin{equation}\lim_{r_d \to\infty} \lim_{t\to\infty} \Xi^{(r_d/t)}(t)
=\operatorname{PPP}\bigl( C Z \sqrt{2}\ee^{-\sqrt{2} x} \,\dd x\bigr) .\vspace*{-3pt}
\end{equation}
\end{prop}

%s3 ###
\section{Proofs} \label{proofs}\vspace*{-3pt}
\mbox{}
\begin{pf*}{Proof of Lemma \protect\ref{comparison}} Theorem \ref
{genealogyabktheorem} describes the genealogies of particles which
fall into compact sets around the level of the maximum but for the
proof of Lemma~\ref{comparison} we need a slight extension in order to
cover the case of sets which are only bounded from below; more
precisely, we claim that for $y\in\R$,
\begin{equation}\label{extension}
\lim_{r_d, r_g \to\infty} \sup_{t > 3 \max\{r_d, r_g \}} \PP[
\exists i,j \in\Sigma_t(y, \infty)\dvtx  Q_{ij}(t) \in(r_d, t-r_g)
] = 0.
\end{equation}
To see this, we recall the following estimate proved by Bramson
\cite{bramson},
Proposition~3:
%e7 ###
%
\begin{equation} \label{boundmax}
\PP\Bigl[\max_{k\leq n(t)} \overline{x_k}(t) \geq Y \Bigr] \leq
\kappa(Y+1)^2\ee^{-\sqrt{2} Y},
\end{equation}
which is valid for $t\geq2, 0< Y < \sqrt{t}$ and $\kappa>0$ a
numerical constant. The bound~\eqref{boundmax} implies in particular that
\begin{equation}\label{newest}
\lim_{Y \to\infty} \sup_{t\geq2} \PP[ \sharp\Sigma_t(Y, \infty)
>0 ] =0.
\end{equation}
For $Y> y$, using the splitting $\Sigma_t(y, \infty)= \Sigma_t(y, Y)
\cup
\Sigma_t(Y, \infty)$, we have the bound
\begin{eqnarray}
&& \PP[ \exists i,j \in\Sigma_t(y, \infty)\dvtx  Q_{ij}(t) \in(r_d,
t-r_g) ]
\nonumber
\\[-9pt]
\\[-9pt]
\nonumber
&&\qquad \leq\PP[ \exists i,j \in\Sigma_t(y, Y)\dvtx  Q_{ij}(t)
\in(r_d, t-r_g) ] + \PP[ \sharp\Sigma_t(Y, \infty) > 0].
\end{eqnarray}
The first term on the right-hand side vanishes, by Theorem \ref
{genealogyabktheorem}, in the limit $t\to\infty$ first and $r_d, r_g
\to\infty$ next, whereas
the second term vanishes, by \eqref{newest}, in the limit $t\to
\infty
$ first and $Y\to\infty$ next: this proves \eqref{extension}.

Let us denote by $A_{t, r_d, r_g}(y, \vare)$ the event
\begin{equation}
\{ \exists i, j \in\Sigma_t(y,\infty)\dvtx  Q_{ij}(t) \in[r_d,
t-r_g] \}.
\end{equation}
By \eqref{extension}, there exists $r_0 = r_0(y, \vare)$ such that, for
$r_d, r_g > r_0$, $\PP[ A_{t, r_d, r_g}^c ] > 1-\vare$. By
definition,
on the event $A_{t, r_d, r_g}^c$, the following equivalence holds for
any $\frac{r_d}{t} \leq q \leq1-\frac{r_g}{t}$:
\begin{equation}
\overline{Q}_{ij}(t) < q \quad\Longleftrightarrow\quad\overline{Q}_{ij}(t) <
\frac{r_d}{t}.
\end{equation}
The assertion of the lemma is now a direct consequence of the
definition of
the thinning $\Xi^{(q)}(t)$ in \eqref{eqnthin}.\vadjust{\goodbreak}
\end{pf*}

To prove Proposition~\ref{simplerpp}, we need some control on
the derivative martingale.

\begin{lem} \label{squarederi}
Let
\begin{equation}
Z^{(2)}(t) \defi\sum_{k\leq n(t)} \bigl\{ \sqrt{2}t - x_k(t) \bigr\}^2
\exp\bigl[- 2 \sqrt{2}\bigl\{ \sqrt{2}t - x_k(t) \bigr\}\bigr].
\end{equation}
For any given $\vare>0$,
\begin{equation}
\lim_{t\to\infty} \PP\bigl[ Z^{(2)}(t) \geq\vare\bigr] = 0.
\end{equation}
\end{lem}

\begin{pf} First, by Bramson's estimate~\cite{bramson}, we may find
$Y = Y(\vare)$ large enough, s.t.
\begin{equation}\label{derivo}
\PP\Bigl[ \max_k x_k(t) - m(t) > Y \Bigr] \leq(1+Y)^2 \ee^{-\sqrt
{2} Y} \leq\vare/2.
\end{equation}
Using this bound, and the Markov inequality, we get
\begin{eqnarray}\label{squareone}\qquad
&& \PP\bigl[ Z^{(2)}(t) \geq\vare\bigr]
\nonumber
\\[-8pt]
\\[-8pt]
\nonumber
&&\qquad\leq\frac{\ee^{t}}{\vare} \E
\bigl[ \bigl\{\sqrt{2}t-x(t)\bigr\}^2 \ee^{-2\sqrt{2} \{\sqrt
{2}t-x(t)\}}; x(t) \leq m(t)+Y \bigr] +\vare/2.
\end{eqnarray}
The first term on the right-hand side is bounded from above by
\begin{eqnarray}\label{squaretwo}
&& \frac{\ee^{t}}{\vare} \int_{( {3}/{(2\sqrt{2})}) \log t-
Y}^{\infty}
x^2 \ee^{-2\sqrt{2} x} \exp\biggl\{-\frac{(\sqrt{2}t -x
)^2}{2t}\biggr\} \frac{\dd x}{\sqrt{2\pi t}} \nonumber\\
&&\qquad \leq\frac{1}{\vare} \int_{ ({3}/{(2\sqrt{2})}) \log t-
Y}^{\infty} x^2 \ee^{-\sqrt{2} x} \ee^{-x^2/2t} \frac{\dd x}{\sqrt
{2\pi
t}} \nonumber\\
&&\qquad \leq\frac{\exp-\sqrt{2}({3}/{(2\sqrt{2})} \log t-
Y)}{\vare} \int_{ ({3}/{(2\sqrt{2})}) \log t- Y}^{\infty} x^2 \ee
^{-x^2/2t} \frac{\dd x}{\sqrt{2\pi t}} \\
&&\qquad \leq\frac{\rho\cdot t^{-3/2}}{ \vare} \int_0^{\infty} x^2 \ee
^{-x^2/2t} \frac{\dd x}{\sqrt{2\pi t}} \nonumber\\
&&\qquad \leq\frac{\rho\cdot t^{-3/2}}{\vare} t \to0 \qquad\mbox{as } t
\to\infty.\nonumber
\end{eqnarray}
This proves the lemma.
\end{pf}

\begin{pf*}{Proof of Proposition \protect\ref{simplerpp}}
We will show the convergence of the Laplace functionals.
For $\phi\dvtx  \R\to\R_+$ measurable with compact support, we claim that
\begin{eqnarray}\label{claimsimpler}
&& \lim_{r_d \to\infty} \lim_{t\to\infty} \E\biggl[ \exp- \int\phi
( x ) \Xi^{(r_d/t)}(t)(\dd x) \biggr]
\nonumber
\\[-8pt]
\\[-8pt]
\nonumber
&&\qquad= \E\biggl[\exp{- C Z
\int\bigl(1-\ee^{-\phi(x)} \bigr) \sqrt{2}\ee^{-\sqrt{2} x} \,\dd
x}\biggr],
\end{eqnarray}
from which the proposition would evidently follow.\vadjust{\goodbreak}

We will prove \eqref{claimsimpler} for simple step functions, that is,
of the form $\phi(x) = \sum_{i=1}^N a_i 1_{A_i}$ for $a_i>0, i=1,
\ldots
, N$ and $A_i = [\underline{A_i}, \overline{A_i}], i=1\ldots N$ disjoint
compact subsets. The extension to the general case of measurable $\phi$
follows by a standard monotone class argument.

We will make use of the splitting
\begin{equation}\label{split}
m(t) = \sqrt{2} r_d + m(t-r_d)+R_t
\end{equation}
for some $R_t = o(1)$ as $t\uparrow\infty$.

We also introduce, for $j \leq n(r_d)$, independent BBMs
$\{x_{k}^{(j)}(t-r_d)\}_{k \leq n_j(t-r_d)}$, and use the abbreviation
\begin{equation}
M_j(t-r_d) \defi\max_{k\leq n_j(t-r_d)} x^{(j)}_{k}(t-r_d) - m(t-r_d).
\end{equation}
Conditionally on everything that happened up to time $r_d$,
the following equality holds in law due to the Markov property and the
definition of the extrema in the $(r_d/t)$-thinning class:
\begin{equation}\label{equalitylaw}
\Xi^{(r_d/t)}(t) \stackrel{(d)}{=} \bigl\{ x_j(r_d) -\sqrt{2} r_d +
M_j(t-r_d)+R_t \bigr\}_{j=1\ldots n(r_d)}.
\end{equation}
Since the $M_j$'s are i.i.d., with $\E_{M(t-r_d)}$ standing for
expectation with respect to $M(t-r_d)$,
\begin{eqnarray}\label{claimsimplerthree}
&&\E\biggl[ \exp- \int\phi( x ) \Xi^{(r_d/t)}(t)(\dd x)
\biggr]
\nonumber
\\[-8pt]
\\[-8pt]
\nonumber
&&\qquad= \E\Biggl[ \prod_{j=1}^{n(r_d)} \E_{M(t-r_d)} \bigl[ \ee^{-\phi
( x_{j}(r_d)-\sqrt{2} r_d + M(t-r_d) +R_t)}\bigr]\Biggr].
\end{eqnarray}
As $t \to\infty$, the variable $M(t-r_d)$ converges weakly to $M$ with
law $\omega$ by~\eqref{travellingone}. Hence
\begin{eqnarray}\label{claimsimplerfour}
&&\lim_{t\to\infty} \E\biggl[ \exp- \int\phi( x ) \Xi
^{(r_d/t)}(t)(\dd x) \biggr]
\nonumber
\\[-8pt]
\\[-8pt]
\nonumber
&&\qquad = \E\Biggl[ \prod_{j=1}^{n(r_d)} \E_{M}
\bigl[ \ee^{-\phi( x_{j}(r_d)-\sqrt{2} r_d+ M )}\bigr]\Biggr].
\end{eqnarray}
Define $y_j(r_d) \defi\sqrt{2} r_d - x_j(r_d)$. We write
\begin{eqnarray}
\E_{M} \bigl[ \ee^{-\phi( -y_j(r_d)+ M )}\bigr] &=& 1 - \E
_{M} \bigl[1- \ee^{-\phi( -y_j(r_d)+ M )}\bigr]
\nonumber
\\[-8pt]
\\[-8pt]
\nonumber
& =:& 1 -
F( -y_j(r_d)),
\end{eqnarray}
and
\begin{equation}
\label{claimsimplerfive}
\eqref{claimsimplerfour}
 = \E\biggl[ \exp\biggl\{ \sum_{j \leq n(r_d)} \log[ 1 - F(
-y_j(r_d)) ] \biggr\} \biggr].
\end{equation}
Observe that
\begin{equation}\min_{j\leq n(r_d)} y_j(r_d) \uparrow\infty\qquad \mbox{a.s.}
\end{equation}
as
$r_d\uparrow\infty$. This implies that
\begin{equation}
\max_{j\leq n(r_d)} F( -y_j(r_d)) \downarrow0.
\end{equation}
Using that $-x -x^2 < \log(1-x) < -x$ for $0<x<1/2$,
for $r_d$ large enough, we
obtain (up to a vanishing error) upper and lower bounds of the form
\begin{eqnarray}\label{claimsimplersix}
 &&\E\biggl[ \exp\biggl\{ - \sum_{j \leq n(r_d)} F( -y_j(r_d))
\biggr\}\biggr]\nonumber\\
&&\qquad \geq\eqref{claimsimplerfive}\\
&&\qquad \geq \E\biggl[ \exp\biggl\{ - \sum_{j \leq n(r_d)} F(
-y_j(r_d)) - F( -y_j(r_d))^2\biggr\} \biggr].\nonumber
\end{eqnarray}
Since $\phi$ is chosen to be a simple step function,
\begin{equation}\label{claimsimplerseven}
F( -y_j(r_d)) = \sum_{i=1}^N (1-\ee^{-a_i}) \int_{A_i +
y_j(r_d)} \,\dd\omega.
\end{equation}
Hence we can make use of the asymptotics \eqref{totheright} to obtain
\begin{eqnarray}\label{claimsimplereight}
 F( -y_j(r_d)) &\sim&
\sum_{i=1}^N (1-\ee^{-a_i}) C \bigl\{ \bigl(\underline{A_i}+y_j(r_d)\bigr) \ee
^{-\sqrt{2}(\underline{A_i}+y_j(r_d)) }
\nonumber
\\[-8pt]
\\[-8pt]
\nonumber
&&\hspace*{48pt}\qquad{} - \bigl(\overline{A_i}+y_j(r_d)\bigr)
\ee
^{-\sqrt{2}(\overline{A_i}+y_j(r_d)) }\bigr\},
\end{eqnarray}
with $\sim$ meaning that the ratio of the left- and right-hand sides
converges to~1, in the limit $r_d \to\infty$, $\PP$-a.s.
We regroup the terms on the right-hand side to get
\begin{eqnarray}\label{claimsimplereight2}
F( -y_j(r_d) ) &\sim& C y_j(r_d) \ee^{-\sqrt{2} y_j(r_d)}
\nonumber
\\[-8pt]
\\[-8pt]
\nonumber
&&{}\times\sum_{i=1}^N (1-\ee^{-a_i}) \{ \ee^{-\sqrt{2} \underline{A_i}} -
\ee^{-\sqrt{2} \overline{A_i}} \} + \mathcal R(y_j(r_d)),
\end{eqnarray}
with $\mathcal R$ containing all the remaining terms; clearly,
\begin{equation}\label{firstrest}
|\mathcal R(y_j(r_d))| \leq\rho\cdot\ee^{-\sqrt{2} y_j(r_d)},
\end{equation}
where $\rho$ depends on the $a_i$ and $A_i$, but not on $y_j(r_d)$. By the
convergence of the derivative martingale as $r_d \uparrow\infty$
[cf. \eqref{eqnderiv}], and the fact that,
in the same limit,
\begin{equation}
\sum_{j\leq n(r_d)} \ee^{-\sqrt{2} y_j(r_d)} \to0,\vadjust{\goodbreak}
\end{equation}
$\PP$-almost surely, we get that
\begin{eqnarray}\label{firstconv}
\lim_{r_d \uparrow\infty} \sum_{j\leq n(r_d)} F\bigl( x_{j}(r_d)-\sqrt
{2} r_d\bigr)&=& C Z \sum_{i=1}^N (1-\ee^{-a_i}) \{ \ee^{-\sqrt{2} \underline
{A_i}} - \ee^{-\sqrt{2} \overline{A_i}} \}
\nonumber\hspace*{-35pt}
\\[-8pt]
\\[-8pt]
\nonumber
& =& C Z \int\bigl( 1- \ee^{-\phi(x)}\bigr) \sqrt{2}\ee^{-\sqrt{2} x}
\,\dd x,\hspace*{-35pt}
\end{eqnarray}
$\PP$-almost surely. This yields the correct asymptotics for the upper bound
in~\ref{claimsimplersix}.

The lower bound in \eqref{claimsimplersix} involves exactly the same
term as the left-hand side of \eqref{firstconv}, and the additional term
\begin{equation}\label{secondconv}
\sum_{j \leq n(r_d)} F\bigl( x_{j}(r_d)-\sqrt{2} r_d\bigr)^2.
\end{equation}
It is straightforward to see that \eqref{secondconv} converges to
zero, as $r_d\uparrow\infty$.
In fact, by the same argument as in \eqref
{claimsimplerseven}--\eqref
{firstrest}, one sees that
\begin{equation}
|\eqref{secondconv}| = O\biggl( \sum_{j \leq n(r_d)} y_j(r_d)^2 \ee
^{-2\sqrt{2} y_j(r_d)}\biggr), \qquad r_d \uparrow\infty.
\end{equation}
With the notation of Lemma~\ref{squarederi},
\begin{equation}
\sum_{j \leq n(r_d)} y_j(r_d)^2 \ee^{-2\sqrt{2} y_j(r_d)} = Z^{(2)}(r_d),
\end{equation}
and this converges to zero in probability, by Lemma~\ref{squarederi}.
Hence, in the limit of large~$r_d$, the lower and upper bounds in
\eqref
{claimsimplersix} coincide, which concludes the proof of the proposition.
\end{pf*}

\begin{pf*}{Proof of Theorem \protect\ref{convergencethinned}} Let $\phi\dvtx
\R
\to\R_+$ be measurable, with compact support.
We need to show that
\begin{eqnarray}
&&\lim_{t\to\infty} \E\biggl[ \exp-\int\phi(x) \Xi
^{(q)}(t)(\dd x)\biggr]
\nonumber
\\[-8pt]
\\[-8pt]
\nonumber
&&\qquad=\E\biggl[\exp{- C Z \int\bigl(1-\ee^{-\phi(x)}
\bigr) \sqrt{2}\ee^{-\sqrt{2} x} \,\dd x}\biggr] ,
\end{eqnarray}
for any $\frac{r_d}{t}\leq q \leq{1-\frac{r_g}{t}}$.
This is straightforward in view of Lemma~\ref{comparison} and
Proposition~\ref{simplerpp} by taking $y$ smaller than the minimum of
the support of $\phi$
and $\vare$ arbitrarily small.
\end{pf*}

%s4 ###
\section{Open problems} \label{secconvection}

%s4.1 ###
\subsection{On the extremal process of BBM}
We consider the following cluster point process.
Let $(\Omega', \mathcal F', P)$, $C>0$ be a probability space, and $Z\dvtx
\Omega'
\to\R_+$
with distribution\vadjust{\goodbreak} as in Theorem~\ref{convergencethinned}.
(Expectation w.r.t.
$P$ will be denoted by
$E$.) Conditionally on a realization of $Z$,
let $(\eta_i; i \in\N)$ be the position of particles generated
according to
a Poisson point process with density
\begin{equation}\label{density}
C Z \bigl( -x \ee^{-\sqrt{2} x} \bigr) \,\dd x
\end{equation}
on the negative axis. For each $i\in\N$, consider independent
branching Brownian motions with drift $-\sqrt{2}$, that is, $\{
x_k^{(i)}(r) -\sqrt{2} r; k \leq n_i(r)\}$, issued on $(\Omega',
\mathcal F', P)$. (``Time'' is denoted here by $r$.)

Remark that for given $i\in\N$,
\begin{equation}\label{driftingoff}
\max_{k \leq n_i(r)} x^{(i)}_k(r)- \sqrt{2} r \to-\infty,
\end{equation}
$P$-almost surely. The branching Brownian motions with drift are then
superimposed on the Poissonian points, that is, the cluster point
process is
given by
\begin{equation}
\qquad\Pi_r \defi\{ \pi_{i, k}(r); i \in\N, k = 1\ldots n_i(r) \}
, \qquad \pi_{i, k}(r) \defi\eta_i + x_k^{(i)}(r)-\sqrt{2} r.
\end{equation}
The existence of the large time limit of $\Pi_r$ is not straightforward.
Due to~\eqref{driftingoff}, only those Poissonian points whose attached
branching Brownian motion performs an unusually large displacement can
contribute to the limiting object. It is thus not clear that one finds any
Poissonian points at all which, together with their cluster of particles,
achieve this feat. The fundamental observation here is that, in virtue of
\eqref{density}, the density
of the Poissonian points on the negative axis grows (slightly faster than)
exponentially when $x\to-\infty$. Together with the work of
Chauvin and Rouault~\cite{chauvinrouault} on branching Brownian motions
conditioned to perform unusually large displacements, this observation
can be exploited to rigorously establish the existence of the point process
$\Pi_r$ in the limit of large times, as well as some of its statistical
properties.
We will report on this in a subsequent paper~\cite{abkthree}.

Here, we only put forward the following conjecture, which appears rather
natural in the light of Theorem~\ref{genealogyabktheorem} and the results
on the paths of extremal particles in BBM established in~\cite{abk}:
\begin{conj} \label{conjabk} In the limit of large times, the
distribution of the extremal process of BBM, $\Xi(t)$ and that of
$\Pi
_r$ coincide, that is,
\begin{equation}\label{equdistr}
\lim_{t\to\infty} \Xi(t) \stackrel{(d)}{=} \lim_{r\to\infty}
\Pi_r .
\end{equation}
In particular, with $\phi\dvtx  \R\to\R_+$ a measurable function with
compact support,
\begin{eqnarray}\label{laplacefull}
&& \lim_{t\to\infty} \E\biggl[\exp\biggl( -\sum_{k\leq n(t)} \phi
\bigl(x_k(t)-m(t)\bigr) \biggr)\biggr]
\nonumber
\\[-8pt]
\\[-8pt]
\nonumber
& &\qquad= \lim_{r\to\infty} E\biggl[ \exp- C Z \int_{-\infty}^0
\bigl( 1 - \ee^{-\psi_r(x)} \bigr) \{- x \ee^{- \sqrt{2} x} \}
\,\dd x \biggr],
\end{eqnarray}
where
\begin{equation}
\ee^{-\psi_r(x)} \defi E\biggl[ \exp\biggl( - \sum_{k\leq n(r)} \phi
\bigl(x + x_k(r) -\sqrt{2} r \bigr) \biggr)\biggr].\vspace*{-3pt}
\end{equation}
\end{conj}

We remark that densities of the form $-x \exp(-\sqrt{2} x) \,\dd x$ on the
negative axis have been conjectured to play an important role in the recent
work by Brunet and Derrida~\cite{derridabrunet}, where the average
size of
the gaps between the $n$th- and ($n+1$)th-leading particle at the edge of
BBM is numerically shown to behave as
\begin{equation}
\frac{1}{n} - \frac{1}{n \log n}+\cdots,
\end{equation}
(which is indeed ``close'' to the average size of the gaps in a PPP
with density $-x \ee^{-\sqrt{2} x} \,\dd x$ on the negative
axis).\vspace*{-3pt}

%s4.2 ###
\subsection{On a conjecture by Lalley and Sellke}
Conjecture~\ref{conjabk} is similar but fundamentally different
from the \textit{Tidal Wave Conjecture} formulated by Lalley and Sellke~\cite{lalleysellke}.
Lalley and Sellke suggested that the Poisson
point process
entering into the construction of $\Pi_r$ should have density
$C Z \ee^{-\sqrt{2} x} \,\dd x$ conditionally on a realization of $Z$
where $C$ is some constant. However, this cannot be correct. We will
show that such a point process
does not exist in the limit $r\to\infty$: the density of the
Poissonian component cannot compensate \eqref{driftingoff} and all the
particles are bound to drift off to $-\infty$. To formulate this
precisely, consider the point process
\begin{equation}\label{lalleysellkepp}
\widetilde\Pi_r\defi\bigl( \tilde\eta_i + x^{(i)}_k(r)- \sqrt{2} r;
i \in\N, k =1,\ldots, n_i(r) \bigr),
\end{equation}
where the $\tilde\eta$'s are points of a PPP with density $C Z \ee
^{-\sqrt{2} x} \,\dd x$, and the $x^{(i)}$'s independent BBMs.\vspace*{-3pt}

\begin{prop} \label{wrongdensity} For given $y\in\R$,
\begin{equation}\label{driftoff}
\lim_{r \to\infty} P\bigl[ \widetilde\Pi_r[y, \infty) \geq1 | Z
\bigr] = 0.\vspace*{-3pt}
\end{equation}
\end{prop}

In order to prove Proposition~\ref{wrongdensity}, we make use of
the following bound established by Bramson:\vspace*{-3pt}
\begin{prop}[(\cite{bramsonmonograph}, Proposition 8.2)] \label
{bramsonremarkablebound}
Let $y_0 < 0$ (strictly). There exists $r_0 = r_0(y_0)$ such that for
$r\geq r_0$, $x \geq m(r)+1$ and $z\defi x-m(r)$,
\begin{equation}\label{uppertight}
P\Bigl[ \max_{k\leq n(r)} x_k(r) \geq x\Bigr] \leq\rho\cdot\ee^r
\int_{y_0}^0 \frac{\ee^{-(x-y)^2/2r}}{\sqrt{2\pi r}} \bigl( 1- \ee^{-2
(y-y_0) z/r}\bigr) \,\dd y,
\end{equation}
where $\rho>0$ is a numerical constant.\vspace*{-3pt}
\end{prop}
Using this with $y_0 := -1$, we obtain the following corollary. (Here
and below, \mbox{$\rho>0$} denotes a numerical constant, not necessarily the
same at different occurrences.)\vadjust{\goodbreak}
\begin{cor}
For $X>1$, and $r\geq r_o = r_o(-1)$,
\begin{eqnarray}\label{uptightzero}
&&P\Bigl[ \max_{k\leq n(t)} x_k(r) -m(r) \geq X \Bigr]
\nonumber
\\[-8pt]
\\[-8pt]
\nonumber
&&\qquad\leq\rho\cdot
X \cdot\exp\biggl(-\sqrt{2} X - \frac{X^2}{2r}+\frac{3 }{2\sqrt{2}}X
\frac{\log r}{r}\biggr).
\end{eqnarray}
\end{cor}

\begin{pf}
According to Proposition~\ref{bramsonremarkablebound}, for $X>1$,
\begin{eqnarray}\label{uptightone}
&&P\Bigl[ \max_{k\leq n(t)} x_k(r) -m(r) \geq X \Bigr]
\nonumber
\\[-8pt]
\\[-8pt]
\nonumber
&&\qquad \leq\rho\cdot
\ee^r \int_{-1}^{0} \frac{\ee^{-(X+m(r)-y)^2/2r}}{\sqrt{2\pi r}} \bigl(
1- \ee^{-2 (y+1) X/r}\bigr) \,\dd y.
\end{eqnarray}
Since $y+1>0$ we have that $1-\ee^{-2 (y+1) X/r} \leq2 (y+1)X/r$.
Using this, the right-hand side of \eqref{uptightone} is at most
\begin{equation}\label{uptightoneone}
\rho\cdot\frac{ X \ee^r}{r} \int_{-1}^{0} (y+1)\frac{\ee
^{-(X+m(r)-y)^2/2r}}{\sqrt{2\pi r}} \,\dd y.
\end{equation}
Expanding the square in the Gaussian density, \eqref{uptightoneone}
is at most
\begin{eqnarray}\label{uptighttwo}
&& \rho\cdot X \cdot\exp\biggl(-\sqrt{2} X - \frac{X^2}{2r}+\frac{3 X
\log r}{2\sqrt{2}r}\biggr)\nonumber\\
&&\quad{}\times \int_{-1}^{0} (y+1) \underbrace{\ee^{Xy/r
+\sqrt{2} y + y({3}/(2\sqrt{2})) (\log r)/r} \ee^{-y^2/2r}}_{\leq1}
\,\dd y
\\
&&\qquad \leq\rho\cdot X \cdot\exp\biggl(-\sqrt{2} X - \frac
{X^2}{2r}+\frac{3 }{2\sqrt{2}}X \frac{\log r}{r}\biggr),\nonumber
\end{eqnarray}
settling the proof of the corollary.
\end{pf}

\begin{pf*}{Proof of Proposition \protect\ref{wrongdensity}}
In view of \eqref
{driftingoff}, it is plain that for any finite set $I\subset\N$
\begin{equation}
\max_{i \in I} \Bigl\{ \tilde\eta_i + \max_{k\leq n_i(r)}
\bigl[x^{(i)}_k(r) -\sqrt{2} r \bigr] \Bigr\} \stackrel{r\uparrow\infty
}{\longrightarrow} -\infty,
\end{equation}
$P$-almost surely. But the number of Poissonian points
$(\tilde\eta_i; i\in\N)$ in the interval $[0, \infty)$ \textit{is}
finite, $P$-almost surely:
this follows from the fact that the density $C Z \ee^{-\sqrt{2}x} \,\dd
x$ is
integrable on $x\in[0, \infty)$. Hence, Proposition \ref
{wrongdensity} will
follow as soon as we prove that
\begin{equation}\label{driftofftwo}
 P\Bigl[ \exists_{i\in\N}\dvtx  \tilde\eta_i + \max_{k\leq n_i(r)}
\bigl\{ x^{(i)}_k(r)- \sqrt{2} r\bigr\} \geq y \mbox{ and } \tilde\eta
_i \in(-\infty, 0) | Z \Bigr] \stackrel{r\uparrow\infty
}{\longrightarrow} 0.\hspace*{-35pt}\vadjust{\goodbreak}
\end{equation}
By the Markov inequality, and using that the BBMs superimposed on the
Poissonian points are identically distributed,
\eqref{driftofftwo} is at most
\begin{equation}\label{driftoffthree}
\int_{-\infty}^0 P\Bigl[ \max_{k\leq n(r)} \bigl\{ x_k(r)- \sqrt{2}
r\bigr\} \geq y - x \Bigr] C Z \ee^{-\sqrt{2} x} \,\dd x.\vspace*{-1pt}
\end{equation}
We rewrite this in terms of $M(r) \defi\max_{k\leq n(r)} \{
x_k(r) - m(r) \} $:
\begin{eqnarray}\label{driftofffour}
\qquad\eqref{driftoffthree} & = &\int_{-\infty}^0 P\biggl[ M(r) \geq y - x +
\frac{3}{2\sqrt{2}} \log r \biggr] C Z \ee^{-\sqrt{2} x} \,\dd x
\nonumber
\\[-10pt]
\\[-10pt]
\nonumber
& = &\bigl(C Z \ee^{-\sqrt{2} y}\bigr) \cdot\frac{1}{r^{3/2}}\int_{y+
{3}/{(2\sqrt{2})} \log r }^\infty P[ M(r) \geq X] \ee
^{\sqrt{2} X} \,\dd X,\vspace*{-1pt}
\end{eqnarray}
the last step by change of variable $y-x+\frac{3}{2\sqrt{2}} \log r
\to X$.

Let us abbreviate $\rho\defi C Z \ee^{-\sqrt{2} y}$. (Note that $y$
and $Z$ are fixed.) For $r$ large enough,
\begin{equation}
y+ \frac{3}{2\sqrt{2}} \log r \geq1,\vspace*{-1pt}
\end{equation}
hence we may use \eqref{uptighttwo} to get that \eqref{driftofffour}
is at most
\begin{eqnarray}\label{driftofffive}\qquad
&& \frac{\rho}{r^{3/2}} \int_{y+ ({3}/{(2\sqrt{2})}) \log r}^\infty X
\exp\biggl( \frac{3 }{2\sqrt{2}}X \frac{\log r}{r}\biggr) \ee^{-
{X^2}/{(2r)}} \,\dd X
\nonumber
\\[-9.5pt]
\\[-9.5pt]
\nonumber
&&\qquad = \frac{\rho}{r^{3/2}} \underbrace{\exp\biggl( \frac{9}{16}
\frac{(\log r)^2}{r} \biggr)}_{= 1+o(1), r\uparrow\infty} \int
_y^\infty\biggl\{ Y+ \frac{3}{2\sqrt{2}} \log r \biggr\} \ee^{-Y^2/2r}
\,\dd Y,\vspace*{-1pt}
\end{eqnarray}
by change of variable $X-\frac{3}{2\sqrt{2}}\log r \to Y$.

It thus remains to control the term
\begin{eqnarray}\label{driftoffsix}\qquad
&& \frac{1}{r^{3/2}} \int_y^\infty\biggl\{ Y+ \frac{3}{2\sqrt{2}} \log r
\biggr\} \ee^{-Y^2/2r} \,\dd Y
\nonumber
\\[-9.5pt]
\\[-9.5pt]
\nonumber
&&\qquad = \frac{1}{r^{3/2}} \int_y^\infty Y \ee^{-Y^2/2r} \,\dd Y+ \frac
{3}{2\sqrt{2}} \cdot\frac{\log r}{r^{3/2}} \int_y^\infty\ee^{-Y^2/2r}
\,\dd Y.\vspace*{-1pt}
\end{eqnarray}
As for the first term on the right-hand side of \eqref{driftoffsix}:
\begin{equation}\label{driftoffseven}\qquad
\frac{1}{r^{3/2}} \int_y^\infty Y \ee^{-Y^2/2r} \,\dd Y = \frac
{1}{\sqrt
{r}} \int_{y/\sqrt{r}}^\infty x \ee^{-x^2/2} \,\dd x \to0 \qquad\mbox{as } r\uparrow\infty.\vspace*{-1pt}
\end{equation}
As for the second term on the right-hand side of \eqref{driftoffsix}:
\begin{eqnarray}\label{driftoffeight}
&&\frac{3}{2\sqrt{2}} \cdot\frac{\log r}{r^{3/2}} \int_y^\infty\ee
^{-Y^2/2r} \,\dd Y
\nonumber
\\[-10pt]
\\[-10pt]
\nonumber
&&\qquad = \frac{3}{2\sqrt{2}} \cdot\frac{\log r}{r} \int
_{y/\sqrt{r}}^\infty\ee^{-x^2/2} \,\dd x \to0\qquad \mbox{as } r\uparrow
\infty.\vspace*{-1pt}
\end{eqnarray}
This proves \eqref{driftofftwo}, settling Proposition \ref
{wrongdensity}.\vadjust{\goodbreak}
\end{pf*}

\begin{rem}
The above computations also suggest that a point process which is
obtained by
superimposing independent BBMs with drift $-\sqrt{2}$ on a PPP with a certain
density exists in the limit of large times if and only if such density is,
up to a (possibly random) constant, $-x \exp(-\sqrt{2} x) \,\dd x$ on the
negative axis.

In fact, a closer look at the above considerations reveals that the
left-hand side of \eqref{driftoffseven} is the leading order of the expected
number of points (of the superimposed point process) which fall into the
subset $[y, \infty)$. Choosing the density of the Poissonian component as
in Conjecture
\ref{conjabk}, \eqref{driftoffseven} would then read $r^{-3/2}
\int
_y^\infty Y^2 \ee^{-Y^2/2r} \,\dd Y$, which indeed remains of order 1 in
the limit $r\to\infty$.
\end{rem}

\section*{Note added in revision} There has been considerable
activity concerning
the extremal process of BBM after this paper was submitted for
publication.
Brunet and Derrida have shown in~\cite{brunetderridatwo} that all
statistical properties of the rightmost points can be extracted from the
traveling wave solutions of the Fisher-KPP equation.
The validity of Conjecture~\ref{conjabk} has been settled in a~paper
of ours~\cite{abkthree}, where it is proved that the extremal process
of branching Brownian motion weakly converges in the limit of large
times to a Poisson cluster process; shortly after that, Aidekon et al.
\cite{aidekonetal}
recovered the same results by means of ``spine techniques.''
The Poissonian structure of the extremal process can also be proved
using the property of \textit{superposability} as observed by Maillard
\cite{maillard}.
This property of the process was conjectured by Brunet and Derrida in
\cite{brunetderridatwo} and proved in~\cite{abkthree}.

\section*{Acknowledgments}
The authors thank \'Eric Brunet and Zakhar Kabluchko for interesting
discussions on the existence of the extremal process of branching
Brownian motion.

% imsref loaded by akundreckaite, 2012-04-03 12:41:12
%

%suskaldyti doi

\printaddresses

\end{document}